\documentclass[12pt,reqno]{amsproc}

\newtheorem{theo}{Theorem}

\newtheorem{lem}{Lemma}

\newcommand\eps\varepsilon
\newcommand\ph\varphi
\newcommand\kap\Lambda

\usepackage{enumerate}
\usepackage{graphicx}
\usepackage{float}
\usepackage{yfonts}
\usepackage{mathrsfs}
\usepackage{amsthm}
\usepackage{amsmath}
\usepackage{amscd}
\usepackage{xcolor}

\begin{document}

\title[On Schauder Bases in Hilbert Space]
{On Schauder Bases in Hilbert Space}

\author[Oleg Zubelevich]{Oleg Zubelevich\\
 Steklov Mathematical Institute of Russian Academy of Sciences \\oezubel@gmail.com
 }

\date{}
\thanks{The research was funded by a grant from the Russian Science
Foundation (Project No. 19-71-30012)}
\subjclass[2000]{46C15, 46B25,47A53, 47B38, 47B48}
\keywords{ Schauder basis, stability of basis, Fredholm operators.}

\begin{abstract}In  this short note we present a far generalization of the following very well-known assertion \cite{halmosh}: assume that we have two orthonormal sequences in a Hilbert space and these sequences are quadratically close to each other. Then if one of these sequences is a basis in the Hilbert space then so is the  other one.

\end{abstract}

\maketitle
\numberwithin{equation}{section}
\newtheorem{theorem}{Theorem}[section]
\newtheorem{lemma}[theorem]{Lemma}
\newtheorem{definition}{Definition}[section]

\section{ The main theorem}Let $H$ be a Hilbert space. To avoid  overloading  the text with details,  we assume that the space $H$ is over the field of reals. But   all the assertions remain valid for the case of
$\mathbb{C}$.

Let $\{e_i\},\quad i\in\mathbb{N}$ be a Schauder basis in $H$ and let $\{e'_i\}$ be its biorthogonal system:
$$(e'_i,e_j)=\delta_{ij}.$$
Here $\delta_{ij}$ stands for the Kronecker delta.

The system $\{e'_i\}$ exists and forms a Schauder basis in $H$ \cite{linden};
$$H\ni x=\sum_{i=1}^\infty(e'_i,x)e_i.$$
Unless otherwise specified,  infinite sums and other limits are regarded in the sense of the strong convergence in $H$.

If $\{e_i\}$ is an orthonormal system then then  $e_i=e'_i$.

\begin{theo}\label{sdgftt} Suppose that a sequence of vectors
\begin{equation}\label{bi}f_1,f_2,\ldots\in H\end{equation} satisfy the following pair of conditions:

1)  the series is convergent:
\begin{equation}\label{zcxvdfg}
\sum_{i,j\in\mathbb{N}}|(e_i,e_j)|\cdot\|f_i-e'_i\|\cdot \|f_j-e'_j\|<\infty;\end{equation}

2) the equality  $\sum_{i=1}^\infty\lambda_if_i=0,\quad \lambda_i\in\mathbb{R}$ implies $$\lambda_i=0,\quad i\in\mathbb{N}.$$ (This property of system (\ref{bi}) is called the $\omega-$independence.)

Then the sequence  (\ref{bi}) is a Schauder basis of $H$.\end{theo}

If sequence (\ref{bi})
is orthonormal  then the condition 2) is satisfied automatically, and if in addition the basis $\{e_i\}$ is orthonormal too then theorem \ref{sdgftt} takes well-known and mentioned in the Abstract form:
the sequence (\ref{bi}) is an orthonormal basis provided
$$\sum_{i\in\mathbb{N}}\|f_i-e_i\|^2<\infty.$$

\subsection{A note on Banach spaces}
Actually the  proof of theorem \ref{sdgftt} (see the next section)  uses the fact that $H$ is a Hilbert space only once, namely in estimate (\ref{sssdg66}).

The main idea of this proof does not  essentially use the Hilbert inner product.

Let $X$ be a reflexive Banach space and let $\{e_i\}$ be its Schauder basis with the corresponding biorthogonal system
\begin{equation}\label{dfgh6y}\{e'_i\}\subset X',\quad (e'_i,e_j)=\delta_{ij}.\end{equation}From \cite{linden} we know that due to reflexivity of $X$ the sequence
(\ref{dfgh6y}) forms a Schauder basis in $X'$.

The proof of the following theorem  is essentially the same as the the proof of theorem
\ref{sdgftt}.
\begin{theo} Let $f_1,f_2,\ldots\in X'$ be an $\omega-$independent sequence.
If
$$\sum_{i=1}^\infty\|e_i\|\cdot\|f_i-e'_i\|<\infty$$ then the sequence $\{f_i\}$ forms a Schauder basis in $X'$.\end{theo}

\section{Proof of theorem \ref{sdgftt}}

Introduce a sequence of finite rank  operators
$$K_n:H\to H,\quad K_nx=\sum_{i=1}^n(f_i-e'_i,x)e_i.$$
For $j>l$ we have an estimate
\begin{equation}\label{sssdg66} \|K_lx-K_jx\|^2\le \|x\|^2\sum_{i,k=l+1,\ldots,j}|(e_i,e_k)|\cdot\|f_i-e'_i\|\cdot \|f_k-e'_k\|.\end{equation}
Thus due to (\ref{zcxvdfg})  the sequence $K_n$ converges to
$$Kx=\sum_{i=1}^\infty(f_i-e'_i,x)e_i$$ in the operator norm.
So that $K$ is a compact operator.

On the other hand we have:
$$K_nx=A_nx-\sum_{i=1}^n(e'_i,x)e_i,\quad A_nx=\sum_{i=1}^n(f_i,x)e_i.$$
Taking the pointwise limit as $n\to\infty$, we obtain a bounded operator $A:H\to H$,
$$Ax=\sum_{i=1}^\infty(f_i,x)e_i=E+K,\quad A_n\to A,$$
where $E$ is the identity map.

Observe that $A$ is a Fredholm operator and its index is equal to zero \cite{Taylor}, \cite{nir}.
This observation implies the following features of the operator $A$.

The kernels of  $A$ and $A^*$ are finite dimensional; their ranges are closed.

The factor spaces  $H/A(H),\quad H/A^*(H)$ are finite dimensional and
\begin{equation}\label{sxdgb677}
\dim\ker A=\dim H/A(H),\quad \dim\ker A^*=\dim H/A^*(H).\end{equation}
From general theory \cite{edv} we also have:
\begin{equation}\label{xsfb678}A(H)=\big(\ker A^*)^\perp.
\end{equation}
By direct calculation one can see that the operators
$$A^*_n:H\to H,\quad A^*_nx=\sum_{i=1}^n(e_i,x)f_i$$  are conjugated to $A_n$.
\begin{lem}\label{sdssgty}
For each $x\in H$ the sequence $A^*_nx$ converges to $A^*x$.\end{lem}
\subsubsection*{Proof of lemma \ref{sdssgty}}

The sequence $\{A^*_nx\}$ is bounded at each point $x$. Indeed, this sequence is weakly convergent:
$$(A^*_nx,y)=(x,A_ny)\to(x,Ay),\quad \forall y\in H.$$

Recall that the sequence $\{e'_i\}$  forms a Schauder basis in $H$.

Let $\mathscr E'=\mathrm{span}\,\{e'_i\}$ denote a linear span of the basis $\{e'_i\}$. The space $\mathscr E'$  is dense in $H$.

For any vector  $x\in \mathscr E'$ the sequence $A^*_nx$ converges, because it stabilizes on each
such  vector.

By the Banach-Steinhaus theorem \cite{edv} the sequence $A^*_n$ converges pointwise to some bounded operator. To finish the proof we pass to the limit in the following equality:
$$(A_n^*x,y)=(x,A_ny).$$


The lemma is proved.

So, $$A^*x=\sum_{i=1}^\infty(e_i,x)f_i.$$
Therefore, by condition 2) of the theorem we have  $\ker A^*=\{0\}.$

Indeed, if $x\in \ker A^*$ then $(e_i,x)=0,\quad \forall i\in\mathbb{N}$. Thus $x\perp \mathscr E=\mathrm{span}\,\{e_i\}.$ But the space $\mathscr E$ is  dense in $H$.

Consequently, on account of equalities  (\ref{sxdgb677}), (\ref{xsfb678})
and by the bounded inverse theorem the operators $A,A^*$ are the linear homeomorphisms of $H$.

Now let us accomplish the proof of the theorem.

Take any vector  $y\in H$.
Then there exists a unique vector $x\in H$ such that
$$y=A^*x=\sum_{i=1}^\infty (e_i,x)f_i,\quad x=(A^*)^{-1}y$$
and
$$y=\sum_{i=1}^\infty \big(e_i,(A^*)^{-1}y\big)f_i.$$
The uniqueness of such an expansion follows from the hypothesis 2) of the theorem.

The theorem is proved.

\end{document}